\documentclass[10pt]{amsart}

\usepackage{graphics}
\usepackage{amsfonts}
\usepackage{amsmath}
\usepackage{amssymb}

\def\RR{{\mathbb R}}
\def\QQ{{\mathbb Q}}
\def\MM{{\mathbb M}}

\def\Z{\mathcal{Z}}
\def\R{\mathcal{R}}
\def\det{{\rm det}}
\def\Aff{{\rm Aff}}
\def\C{C}

\def\dto{{\scriptstyle\buildrel{\scriptstyle\longrightarrow}
\over{{}_{\scriptstyle\longrightarrow}}}}
\def\eps{\varepsilon}
\def\E{{\mathcal E}}
\def\I{{\mathcal E}}
\def\inff{\mathop{\inf\limits_{a\in\RR^2}}}
\def\infff{\mathop{\inff\limits_{b\in\RR^3}}}

\newtheorem{theorem}{Theorem}[section]
\newtheorem{lemma}[theorem]{Lemma}
\newtheorem{proposition}[theorem]{Proposition}
\newtheorem{corollary}[theorem]{Corollary}

\theoremstyle{definition}
\newtheorem{definition}[theorem]{Definition}

\theoremstyle{remark}
\newtheorem{remark}[theorem]{Remark}

\title[The nonlinear membrane energy]{The nonlinear membrane energy: variational derivation under the constraint \boldmath$``\det\nabla u>0"$\unboldmath}

\author{Omar Anza Hafsa}
\address{Institut f\"ur Mathematik, Universit\"at Z\"urich, Winterthurerstrasse 190, CH-8057 Z\"urich, Switzerland.}
\email{anza@math.unizh.ch}

\author{Jean-Philippe Mandallena}
\address{``Equipe AVA (Analyse Variationnelle et Applications)", Centre Universitaire de Formation et de Recherche de N\^\i mes, Site des Carmes, Place Gabriel P\'eri - Cedex 01 - 30021 N\^\i mes, France.\newline \indent I3M (Institut de Math\'ematiques et Mod\'elisation de Montpellier) UMR - CNRS 5149, Universit\'e Montpellier II, Place Eug\`ene Bataillon, 34090 Montpellier, France.}
\email{jean-philippe.mandallena@unimes.fr}

\begin{document}

\begin{abstract}
In \cite{anzman-det} we gave a variational definition of the nonlinear membrane energy under the constraint ``$\det\nabla u\not=0$". In this paper we obtain the nonlinear membrane energy under the more realistic constraint ``$\det\nabla u>0$". 
\end{abstract}

\maketitle


\section{Introduction}

Consider an elastic material occupying in a reference configuration the bounded open set $\Sigma_\eps\subset\RR^3$ given by
$$
\Sigma_\eps:=\Sigma\times\left]-{\eps\over 2},{\eps\over 2}\right[,
$$
where $\eps>0$ is very small and $\Sigma\subset\RR^2$ is Lipschitz, open and bounded. A point of $\Sigma_\eps$ is denoted by $(x,x_3)$ with $x\in\Sigma$ and $x_3\in]-{\eps\over 2},{\eps\over 2}[$.
Let
$
W:\MM^{3\times 3}\to[0,+\infty]
$ 
be the stored-energy function supposed to be {\em continuous} and {\em coercive}, i.e., $W(F)\geq C|F|^p$ for all $F\in\MM^{3\times 3}$ and some $C>0$. In order to take into account the important physical properties that the interpenetration of matter does not occur and that an infinite amount of energy is required to compress a finite volume into zero volume, i.e., 
$$
W(F)\to+\infty\ \hbox{ as }\ \det F\to 0,
$$
where $\det F$ denotes the determinant of the $3\times 3$ matrix $F$, we assume that:
\begin{equation}\label{Nim}
\hbox{\em $W(F)=+\infty$ if and only if $\det F\leq 0$};
\end{equation}
\begin{equation}\label{Icv}
\hbox{\em for every $\delta>0$, there exists $c_\delta>0$ such that for all $F\in\MM^{3\times 3}$},
\end{equation}
$$
\hbox{\em if }\det F\geq\delta\hbox{\em\ then } W(F)\leq c_\delta(1+|F|^p).
$$
Our goal is to show that as $\eps\to 0$ the three-dimensional free energy functional ${E}_\eps:W^{1,p}(\Sigma_\eps;\RR^3)\to[0,+\infty]$ (with $p>1$) defined by
\begin{equation}\label{SEF}
{E}_\eps(u):={1\over\eps}\int_{\Sigma_\eps}W\big(\nabla u(x,x_3)\big)dxdx_3
\end{equation}
converges in a variational sense (see Definition \ref{variationalconvergence}) to the two-dimensional free energy functional ${E}_{\rm mem}:W^{1,p}(\Sigma;\RR^3)\to[0,+\infty]$ given by
\begin{equation}\label{MSEF}
{E}_{\rm mem}(v):=\int_\Sigma W_{\rm mem}\big(\nabla v(x)\big)dx
\end{equation}
with $W_{\rm mem}:\MM^{3\times 2}\to[0,+\infty]$. Usually, $E_{\rm mem}$ is called the nonlinear membrane energy associated with the two-dimensional elastic material with respect to the reference configuration $\Sigma$. Furthermore we wish to give a representation formula for $W_{\rm mem}$. 

To our knowledge, the problem of giving a variational definition of the nonlinear membrane energy was studied for the first time by Percivale in \cite{percivale}. His paper deals with the constraint ``$\det\nabla u>0$" but seems to contain some mistakes (it never was published). Nevertheless, Percivale introduced the ``good" formula for $W_{\rm mem}$, i.e., $W_{\rm mem}=\mathcal{Q}W_0$ where $W_0$ is given by (\ref{PercivaleFormula}) and $\mathcal{Q}W_0$ denotes the quasiconvex envelope of $W_0$.  
Then, in \cite{ledretraoult1} Le Dret and Raoult gave a complete proof of percivale's conjecture in the simpler case where $W$ is of $p$-polynomial growth, i.e., $W(F)\leq c(1+|F|^p)$ for all $F\in\MM^{3\times 3}$ and some $c>0$. Although the $p$-polynomial growth case is not compatible with (\ref{Nim}) and (\ref{Icv}) their paper established a suitable framework to deal with dimensional reduction problems (it is the point of departure of many works on the subject). After Percivale, Ben Belgacem also considered the constraint ``$\det\nabla u>0$". In \cite[Theorem 1]{benbelgacem1} he announced to have succeed to handle (\ref{Nim}) and (\ref{Icv}). In \cite{benbelgacem}, which is the paper corresponding to the note \cite{benbelgacem1}, the statement \cite[Theorem 1]{benbelgacem1} is partly proved (however, a more detailled proof, but not fully complete, can be found in his thesis \cite{benbelgacem-thesis}). Moreover, for Ben Belgacem $W_{\rm mem}=\mathcal{Q}\mathcal{R}W_0$ where $\mathcal{R}W_0$ denotes the rank one convex envelope of $W_0$ (in fact, as we proved in \cite{anzman2,anzman-det}, $\mathcal{Q}\mathcal{R}W_0=\mathcal{Q}W_0$). Nevertheless, Ben Belgacem's thesis highlighted the role of approximation theorems for Sobolev functions by smooth immersions in the studying of the passage 3D-2D in presence of (\ref{Nim}) and (\ref{Icv}).  Recently, in \cite{anzman-det} we gave a variational definition of the nonlinear membrane energy under the constraint ``$\det\nabla u\not=0$". In the present paper, using the same method as in \cite{anzman-det} and some results of Ben Belgacem's thesis (mainly, Theorem \ref{LemmaBBB} and Lemma \ref{LeMMaBiS}), we obtain the nonlinear membrane energy under the more realistic constraint ``$\det\nabla u>0$".

\smallskip

An outline of the paper is as follows. The variational convergence of $E_\eps$ to $E_{\rm mem}$ as $\eps\to 0$ as well as a representation formula for $W_{\rm mem}$ are given by Corollary \ref{corollary} in Sect. 2.4. Corollary \ref{corollary} is a consequence of Theorems \ref{first_main_result}, \ref{AdDitioNal} and \ref{anzman}. Roughly,  Theorems \ref{first_main_result} and \ref{AdDitioNal} establish the existence of the variational limit of ${E}_\eps$ as $\eps\to 0$ (see Sect. 2.2), and Theorem \ref{anzman} gives an integral representation for the corresponding variational limit, and so a representation formula for $W_{\rm mem}$ (see Sect. 2.3). In fact, Theorem \ref{anzman}  is obtained from Theorem \ref{AdDitioNal} which furnishes a ``simplified" formula for the variational limit.

Theorem \ref{first_main_result} is proved in Section 4. The principal ingredients are Theorem \ref{AdDitioNal} and Theorem  \ref{basictheorem} whose proof (given in Section 3) uses an interchange theorem of infimum and integral that we obtained in \cite{anzman}. 
(Note that the techniques used to prove Theorems \ref{first_main_result} and \ref{basictheorem} are the same as in \cite[Sections 3 and 4]{anzman-det}.) 

Theorem \ref{AdDitioNal}  is proved is Section 5. The main arguments are two approximation theorems developed by Ben Belgacem-Bennequin (see \cite{benbelgacem-thesis}) and Gromov-Eliashberg (see \cite{gromovEliashberg}). These theorems are stated in Appendix A. 

Theorem \ref{anzman} is proved in \cite[Appendix A]{anzman-det} (see also \cite{anzman2}).


\section{Results}


\subsection{Variational convergence} 
To accomplish our asymptotic analysis, we use the notion of convergence introduced by Anzellotti, Baldo and Percivale in \cite{anzebalper} in order to deal with dimension reduction problems in mechanics. Let $\pi=\{\pi_\eps\}_\eps$ be the family of maps $\pi_\eps:W^{1,p}(\Sigma_\eps;\RR^3)\to W^{1,p}(\Sigma;\RR^3)$ defined by
$$
\pi_\eps(u):={1\over\eps}\int_{-{\eps\over 2}}^{\eps\over 2}u(\cdot,x_3)dx_3.
$$
\begin{definition}\label{variationalconvergence}
We say that ${E}_\eps$ $\Gamma(\pi)$-converges to ${E}_{\rm mem}$ as $\eps\to 0$, and we write ${E}_{\rm mem}=\Gamma(\pi)\hbox{\rm -}\lim_{\eps\to 0}{E}_\eps$, if the following two assertions hold{\rm:}
\begin{itemize}
\item[(i)] for all $v\in W^{1,p}(\Sigma;\RR^3)$ and all $\{u_\eps\}_\eps\subset W^{1,p}(\Sigma_\eps;\RR^3)$, 
$$
\hbox{if }\pi_\eps(u_\eps)\to v\hbox{ in }L^{p}(\Sigma;\RR^3)\hbox{ then }{E}_{\rm mem}(v)\leq\liminf_{\eps\to 0}{E}_\eps(u_\eps);
$$
\item[(ii)] for all $v\in W^{1,p}(\Sigma;\RR^3)$, there exists $\{u_\eps\}_\eps\subset W^{1,p}(\Sigma_\eps;\RR^3)$ such that{\rm:} 
$$
\pi_\eps(u_\eps)\to v\hbox{ in }L^{p}(\Sigma;\RR^3)\hbox{ and }{E}_{\rm mem}(v)\geq\limsup_{\eps\to 0}{E}_\eps(u_\eps).
$$
\end{itemize}
\end{definition}
\noindent In fact, Definition \ref{variationalconvergence} is a variant of De Giorgi's $\Gamma$-convergence. This is made clear by Lemma \ref{link}. Consider $\mathcal{\I}_\eps:W^{1,p}(\Sigma;\RR^3)\to[0,+\infty]$ defined by
$$
\mathcal{\I}_\eps(v):=\inf\Big\{{E}_\eps(u):\pi_\eps(u)=v\Big\}.
$$
\begin{definition}\label{gammavariationalconvergence}
We say that $\mathcal{\I}_\eps$ $\Gamma$-converges to ${E}_{\rm mem}$ as $\eps\to 0$, and we write ${E}_{\rm mem}=\Gamma\hbox{\rm -}\lim_{\eps\to 0}\mathcal{\I}_\eps$ if for every $v\in W^{1,p}(\Sigma;\RR^3)$,
$$
\left(\Gamma\hbox{-}\liminf_{\eps\to 0}\mathcal{\I}_\eps\right)(v)=\left(\Gamma\hbox{-}\limsup_{\eps\to 0}\mathcal{\I}_\eps\right)(v)=E_{\rm mem}(v),
$$
where  
$
\left(\Gamma\hbox{-}\liminf_{\eps\to 0}\mathcal{\I}_\eps\right)(v):=\inf\big\{\liminf_{\eps\to 0}\mathcal{\I}_\eps(v_\eps):v_\eps\to v\hbox{ in }L^{p}(\Sigma;\RR^3)\big\}
$
and
$
\left(\Gamma\hbox{-}\limsup_{\eps\to 0}\mathcal{\I}_\eps\right)(v):=\inf\big\{\limsup_{\eps\to 0}\mathcal{\I}_\eps(v_\eps):v_\eps\to v\hbox{ in }L^{p}(\Sigma;\RR^3)\big\}.
$
\end{definition}
For a deeper discussion of the $\Gamma$-convergence theory we refer to the book \cite{dalmaso0}. Clearly, Definition \ref{gammavariationalconvergence} is equivalent to assertions (i) and (ii) in definition \ref{variationalconvergence} with ``$\pi(u_\eps)\to v$'' replaced by ``$v_\eps\to v$''. It is then obvious that
\begin{lemma}\label{link}
${E}_{\rm mem}=\Gamma(\pi)\hbox{\rm -}\lim_{\eps\to 0}{E}_\eps$ if and only if ${E}_{\rm mem}=\Gamma\hbox{\rm -}\lim_{\eps\to 0}\mathcal{\I}_\eps$.
\end{lemma}

                                    
The $\Gamma(\pi)$-convergence of ${E}_\eps$ in (\ref{SEF}) to $E_{\rm mem}$ in (\ref{MSEF}) as $\eps\to 0$ as well as a representation formula for $W_{\rm mem}$ are given by Corollary \ref{corollary}. It is a consequence of Theorems \ref{first_main_result}, \ref{AdDitioNal} and \ref{anzman}. Roughly,  Theorems \ref{first_main_result} and \ref{AdDitioNal} establish the existence of the $\Gamma(\pi)$-limit of ${E}_\eps$ as $\eps\to 0$ (see Sect. 2.2), and Theorem \ref{anzman} gives an integral representation for the corresponding $\Gamma(\pi)$-limit, and so a representation formula for $W_{\rm mem}$ (see Sect. 2.3). 


\subsection{\boldmath$\Gamma$\unboldmath-convergence of \boldmath$\mathcal{\I}_\eps$\unboldmath\ as \boldmath$\eps\to 0$\unboldmath} 
Denote by $\C^1(\overline{\Sigma};\RR^3)$ the space of all restrictions to $\overline{\Sigma}$ of $C^1$-differentiable functions from $\RR^2$ to $\RR^3$, and set 
$$
\C^1_*(\overline{\Sigma};\RR^3):=\Big\{v\in\C^1(\overline{\Sigma};\RR^3):\partial_1v(x)\land\partial_2 v(x)\not=0\hbox{ for all }x\in\overline{\Sigma}\Big\},
$$
where $\partial_1v(x)$ (resp. $\partial_2 v(x)$) denotes the partial derivative of $v$ at $x=(x_1,x_2)$ with respect to $x_1$ (resp. $x_2$). (In fact, $\C^1_*(\overline{\Sigma};\RR^3)$ is the set of all $C^1$-immersions from $\overline{\Sigma}$ to $\RR^3$.) Let ${\mathcal E}:W^{1,p}(\Sigma;\RR^3)\to[0,+\infty]$ be defined by
$$
{\mathcal E}(v):=\left\{
\begin{array}{cl}
\displaystyle\int_\Sigma W_0\big(\nabla v(x)\big)dx&\hbox{if }v\in\C^1_*(\overline{\Sigma};\RR^3)\\
+\infty&\hbox{otherwise,}
\end{array}
\right.
$$
where $W_0:\MM^{3\times 2}\to[0,+\infty]$ is given by
\begin{equation}\label{PercivaleFormula}
W_0(\xi):=\inf_{\zeta\in\RR^3}W(\xi\mid\zeta)
\end{equation}
with $(\xi\mid\zeta)$ denoting the element of $\MM^{3\times 3}$ corresponding to $(\xi,\zeta)\in\MM^{3\times 2}\times\RR^3$. (As $W$ is coercive, it is easy to see that {\em $W_0$ is coercive}, i.e., $W_0(\xi)\geq C|\xi|^p$ for all $\xi\in\MM^{3\times 2}$ and some $C>0$.) The following lemma gives three elementary properties of $W_0$ (the proof is left to the reader). Note that conditions (\ref{Nim}) and (\ref{Icv}) imply $W_0$ is not of $p$-polynomial growth. 
\begin{lemma}\label{propertiesofW_0} 
Denote by $\xi_1\land\xi_2$ the cross product of vectors $\xi_1,\xi_2\in\RR^3$.
\begin{itemize}
\item[(i)] $W_0$ is continuous.
\item[(ii)] If {\rm(\ref{Nim})} holds then \hbox{$W_0(\xi_1\mid\xi_2)=+\infty$ if and only if $\xi_1\land\xi_2=0$.}
\item[(iii)] If {\rm(\ref{Icv})} holds then{\rm:} 
\begin{equation}\label{Icvbis}
\hbox{for all $\delta>0$, there exists $c_\delta>0$ such that for all $\xi=(\xi_1\mid\xi_2)\in\MM^{3\times 2}$},
\end{equation}
$$
\hbox{if }|\xi_1\land\xi_2|\geq\delta\hbox{ then } W_0(\xi)\leq c_\delta(1+|\xi|^p). 
$$
\end{itemize}
\end{lemma}
Taking Lemma \ref{link} into account, we see that the existence of the $\Gamma(\pi)$-limit of $E_\eps$ as $\eps\to 0$ follows from Theorem \ref{first_main_result}.
\begin{theorem}\label{first_main_result}
Let assumptions {\rm(\ref{Nim})} and  {\rm(\ref{Icv})} hold. Then $\Gamma\hbox{\rm -}\lim_{\eps\to 0}\mathcal{\I}_\eps=\overline{{\mathcal E}}$ with  $\overline{{\mathcal E}}:W^{1,p}(\Sigma;\RR^3)\to[0,+\infty]$ given by
$$
\overline{{\mathcal E}}(v):=\inf\left\{\liminf_{n\to+\infty}{\mathcal E}(v_n):W^{1,p}(\Sigma;\RR^3)\ni v_n\to v\hbox{ in }L^{p}(\Sigma;\RR^3)\right\}.
$$
\end{theorem}
The proof of Theorem \ref{first_main_result} is established in Section 4. It uses Theorem  \ref{basictheorem} (see Section 3) and Theorem \ref{AdDitioNal}.
\begin{theorem}\label{AdDitioNal}
If {\rm(\ref{Icvbis})} holds then $\overline{{\mathcal E}}(v)={\mathcal I}(v)$ for all $v\in W^{1,p}(\Sigma;\RR^3)$, where ${\mathcal I}:W^{1,p}(\Sigma;\RR^3)\to[0,+\infty]$ is given by
$$
{\mathcal I}(v):=\inf\left\{\liminf_{n\to+\infty}\int_\Sigma W_0\big(\nabla v_n(x)\big)dx:W^{1,p}(\Sigma;\RR^3)\ni v_n\to v\hbox{ in }L^{p}(\Sigma;\RR^3)\right\}.
$$
\end{theorem}
Theorem \ref{AdDitioNal} is proved in Section 6 by using two approximation theorems developed by Ben Belgacem-Bennequin (see \cite{benbelgacem-thesis}) and Gromov-Eliashberg (see \cite{gromovEliashberg}). These theorems are stated in Appendix A. 

                                                     
\subsection{Integral representation of \boldmath$\mathcal{I}$\unboldmath} From now on, given a bounded open set $D\subset\RR^2$ with $|\partial D|=0$, we denote by $\Aff(D;\RR^3)$ the space of all continuous piecewise affine functions from $D$ to $\RR^3$, i.e., {\em $v\in\Aff(D;\RR^3)$ if and only if $v$ is continuous and there exists a finite family $(D_i)_{i\in I}$ of open disjoint subsets of $D$ such that $|\partial D_i|=0$ for all $i\in I$, $|D\setminus \cup_{i\in I} D_i|=0$ and for every $i\in I$, $\nabla v(x)=\xi_i$ in $D_i$ with $\xi_i\in\MM^{3\times 2}$} (where $|\cdot|$ denotes the Lebesgue measure in $\RR^2$). Define $\Z W_0:\MM^{3\times 2}\to[0,+\infty]$ by 
\begin{equation}\label{DefinitionofZW0}
\Z W_0(\xi):=\inf\left\{\int_Y W_0\big(\xi+\nabla\phi(y)\big)dy:\phi\in \Aff_0(Y;\RR^3)\right\}
\end{equation}
where $Y:=]0,1[^2$ and $\Aff_0(Y;\RR^3):=\{\phi\in \Aff(Y;\RR^3):\phi=0\hbox{ on }\partial Y\}$. (As $W_0$ is coercive, it is easy to see that {\em $\Z W_0$ is coercive}.) Recall the definitions of quasiconvexity and quasiconvex envelope:
\begin{definition}\label{DEFofQUasiRankoNEConvexityandEnvelOPe}
Let $f:\MM^{3\times 2}\to[0,+\infty]$ be a Borel measurable function.
\begin{itemize}
\item[(i)] We say that $f$ is quasiconvex if for every $\xi\in\MM^{3\times 2}$, every bounded open set $D\subset\RR^2$ with $|\partial D|=0$ and every $\phi\in W^{1,\infty}_0(D;\RR^3)$,
$$
f(\xi)\leq{1\over|D|}\int_D f(\xi+\nabla\phi(x))dx.
$$ 
\item[(ii)] By the quasiconvex envelope of $f$, we mean the unique function (when it exists) $\mathcal{Q}f:\MM^{3\times 2}\to[0,+\infty]$ such that:
\begin{itemize}
\item[-] $\mathcal{Q}f$ is Borel measurable, quasiconvex and $\mathcal{Q}f\leq f$;
\item[-] for all $g:\MM^{3\times 2}\to[0,+\infty]$, if $g$ is Borel measurable, quasiconvex and $g\leq f$, then $g\leq\mathcal{Q}f$.
\end{itemize} 
(Usually, for simplicity, we say that $\mathcal{Q}f$ is the greatest quasiconvex function which less than or equal to $f$.) 
\end{itemize}
\end{definition}
Under  {\rm(\ref{Icvbis})}, we proved that  $\Z W_0$ is of $p$-polynomial growth and so continuous (see \cite[Propositions A.3 and A.1(iii)]{anzman-det}) and that $\Z W_0$ is the quasiconvex envelope of $W_0$, i.e., $\Z W_0=\mathcal{Q}W_0$ (see \cite[Proposition A.5]{anzman-det}). Taking Theorems \ref{first_main_result} and \ref{AdDitioNal} together with Lemmas \ref{link} and \ref{propertiesofW_0}(iii) into account, we see that Theorem \ref{anzman} gives an integral representation for the $\Gamma(\pi)$-limit of $E_\eps$ as $\eps\to 0$ as well as a representation formula for $W_{\rm mem}$.
\begin{theorem}\label{anzman}
If {\rm(\ref{Icvbis})} holds then for every $v\in W^{1,p}(\Sigma;\RR^3)$,
$$
{\mathcal I}(v)=\int_\Sigma \mathcal{Q}W_0\big(\nabla v(x)\big)dx.
$$
\end{theorem}
 Theorem \ref{anzman}  is proved in \cite[Appendix A]{anzman-det} (see also \cite{anzman2}).
 
\subsection{\boldmath$\Gamma(\pi)$\unboldmath-convergence of \boldmath${E}_\eps$\unboldmath\ to \boldmath${E}_{\rm mem}$\unboldmath\ as \boldmath$\eps\to 0$\unboldmath} 
According to Lemmas \ref{link} and Lemma \ref{propertiesofW_0}(iii), a direct consequence of Theorems \ref{first_main_result}, \ref{AdDitioNal} and \ref{anzman} is the following.
\begin{corollary}\label{corollary}
Let assumptions {\rm(\ref{Nim})} and {\rm(\ref{Icv})} hold. Then as $\eps\to 0$, ${E}_\eps$ in {\rm(\ref{SEF})} $\Gamma(\pi)$-converge to ${E}_{\rm mem}$ in {\rm(\ref{MSEF})} with $W_{\rm mem}=\mathcal{Q}W_0$. 
\end{corollary}
\begin{remark}
Corollary \ref{corollary} can be applied when $W:\MM^{3\times 3}\to[0,+\infty]$ is given by
$$
W(F):=
 h(\det F)+|F|^p,
$$
where $h:\RR\to[0,+\infty]$  is a  continuous function such that:
\begin{itemize}
\item[-] {$h(t)=+\infty$ if and only if $t\leq0$};
\item[-] {for every $\delta>0$, there exists $r_\delta>0$ such that $h(t)\leq r_\delta$ for all $t\geq\delta$}.
\end{itemize}
\end{remark}


\section{Representation of ${\mathcal E}$} 

The goal of this section is to show Theorem  \ref{basictheorem}. To this end, we begin by proving two lemmas. 

For every $v\in\C^1_*(\overline{\Sigma};\RR^3)$ and $j\geq 1$, we define the multifunction $\Lambda^j_v:\overline{\Sigma}\dto\RR^3$ by
$$
\Lambda_v^j(x):=\left\{\zeta\in\RR^3: \det(\nabla v(x)\mid\zeta)\ge {1\over j}\right\}.
$$
 \begin{lemma}\label{Lemma2}
Let $v\in\C^1_*(\overline{\Sigma};\RR^3)$. Then{\rm:}
\begin{itemize}
\item[{\rm(i)}] for every $j\geq 1$, $\Lambda_v^j$  is a nonempty convex closed-valued lower semicontinuous\footnote{A multifunction $\Lambda:\overline{\Sigma}\to\RR^3$ is said to be lower semicontinuous if for every closed subset $X$ of $\RR^3$, every $x\in\overline{\Sigma}$ and every $\{x_n\}_{n\geq 1}\subset\overline{\Sigma}$ such that $|x_n-x|\to 0$ as $n\to+\infty$ and $\Lambda(x_n)\subset X$ for all $n\geq 1$, we have $\Lambda(x)\subset X$ (see \cite{aubinfrankowska} for more details).} multifunction{\rm;}
\item[(ii)] for every $x\in\overline{\Sigma}$, $\Lambda_v^1(x)\subset\cdots\subset\Lambda_v^j(x)\subset\cdots\subset\cup_{j\ge 1}\Lambda_v^j(x)=\Lambda_v(x)$, where $\Lambda_v(x):=\{\zeta\in\RR^3: \det(\nabla v(x)\mid\zeta)>0\}$.
\end{itemize}
\end{lemma}
\begin{proof} (ii) is obvious. Prove then (i). Let $j\geq 1$. It is easy to see that for every $x\in\overline{\Sigma}$, $\Lambda_v^j(x)$ is nonempty, convex and closed. Let $X$ be a closed subset of $\RR^3$, let $x\in\overline{\Sigma}$, and let $\{x_n\}_{n\geq 1}\subset\overline{\Sigma}$ such that $|x_n-x|\to 0$ as $n\to+\infty$ and $\Lambda^j_{v}(x_n)\subset X$ for all $n\geq 1$. Let $\zeta\in \Lambda_v^j(x)$ and let $\{\zeta_m\}_{m\geq 1}\subset\RR^3$ be given by $\zeta_m:=\zeta+{1\over m}\zeta$. Then, for every $m\geq 1$,
\begin{equation}\label{DetEq}
\det\big(\nabla v(x)\mid\zeta_m\big)=\det\big(\nabla v(x)\mid\zeta\big)+{1\over m}\det\big(\nabla v(x)\mid\zeta\big)\geq{1\over j}+{1\over mj}.
\end{equation}
Fix any $m\geq 1$. Since $\det(\nabla v(x_n)\mid\zeta_m)\to\det(\nabla v(x)\mid\zeta_m)$ as $n\to +\infty$, using (\ref{DetEq}) we see that $\det(\nabla v(x_{n_0})\mid\zeta_m)>{1\over j}$ for some $n_0\geq 1$, so that $\zeta_m\in\Lambda_v^j(x_{n_0})$. Thus $\zeta_m\in X$ for all $m\geq 1$. As $X$ is closed we have $\zeta=\lim_{m\to +\infty}\zeta_m\in X$.
\end{proof}

In the sequel, given $\Lambda:\overline{\Sigma}\dto\RR^3$ we set
$$
C(\overline{\Sigma};\Lambda):=\Big\{\phi\in C\big(\overline{\Sigma};\RR^3\big):\phi(x)\in\Lambda(x)\hbox{ for all }x\in\overline{\Sigma}\Big\},
$$ 
where $C(\overline{\Sigma};\RR^3)$ denotes the space of all continuous functions from $\overline{\Sigma}$ to $\RR^3$. 
\begin{lemma}\label{Lemma3}
Given $v\in\C^1_*(\Sigma;\RR^3)$ and $j\geq 1$, if {\rm(\ref{Icv})} holds, then
$$
\inf_{\phi\in C(\overline{\Sigma};\Lambda_v^j)}\int_{\Sigma}W\big(\nabla v(x)\mid\phi(x)\big)dx=\int_\Sigma \inf_{\zeta\in\Lambda_v^j(x)}W\big(\nabla v(x)\mid\zeta\big)dx.
$$
\end{lemma} 
To prove Lemma \ref{Lemma3} we need the following interchange theorem of infimum and integral (that we proved in \cite[Corollary 5.4]{anzman}).
\begin{theorem}\label{AHM} Let $\Gamma:\overline{\Sigma}\dto\RR^3$ and let $f:\overline{\Sigma}\times\RR^3\to[0,+\infty]$. Assume that{\rm:} 
\begin{itemize}
\item[(H$_1$)]  $f$ is a Carath\'eodory integrand{\rm;}
\item[(H$_2$)] $\Gamma$ is a nonempty convex closed-valued lower semicontinuous multifunction{\rm;}
\item[(H$_3$)] $C(\overline{\Sigma};\Gamma)\not=\emptyset$ and for every $\phi,\hat\phi\in C(\overline{\Sigma};\Gamma)$,
$$
\int_\Sigma \max_{\alpha\in[0,1]}f\big(x,\alpha\phi(x)+(1-\alpha)\hat\phi(x)\big)dx<+\infty.
$$
\end{itemize}
Then,
$$
\inf_{\phi\in C(\overline{\Sigma};\Gamma)}\int_{\Sigma}f\big(x,\phi(x)\big)dx=\int_\Sigma\inf_{\zeta\in\Gamma(x)}f(x,\zeta)dx.
$$
\end{theorem}

\smallskip

\noindent{\em Proof of Lemma {\rm\ref{Lemma3}.} }Since $W$ is continuous, (H$_1$) holds with $f(x,\zeta)=W(\nabla v(x)\mid\zeta)$. Lemma \ref{Lemma2} shows that (H$_2$) is satisfied with $\Gamma=\Lambda_v^j$, and $C(\overline{\Sigma};\Lambda^j_v)\not=\emptyset$ (for example $\Phi:\overline{\Sigma}\to\RR^3$ defined by (\ref{ExampleOfContinuousSelection}) belongs to $C(\overline{\Sigma};\Lambda^j_v)$). Given $\phi,\hat\phi\in C(\overline{\Sigma};\Lambda^j_v)$, it is clear that 
$
\det(\nabla v(x)\mid\alpha\phi(x)+(1-\alpha)\hat\phi(x))\geq{1/ j}
$
for all $\alpha\in[0,1]$  and all $x\in\overline{\Sigma}$. By (\ref{Icv}) there exists $c>0$ depending only on $j$, $v$, $\phi$ and $\hat\phi$ such that
$
W(\nabla v(x)\mid\alpha\phi(x)+(1-\alpha)\hat\phi(x))\leq c
$
for all $x\in\overline{\Sigma}$. Thus (H$_3$) is verified with $f(x,\zeta)=W(\nabla v(x)\mid\zeta)$ and $\Gamma=\Lambda_v^j$, and Lemma \ref{Lemma3} follows from Lemma \ref{AHM}.\hfill$\square$
  
 \medskip              

Here is our (non integral) representation theorem for ${\mathcal E}$.
\begin{theorem}\label{basictheorem} If {\rm(\ref{Nim})} and  {\rm(\ref{Icv})} hold, then for every $v\in C^1_*(\overline{\Sigma};\RR^3)$,
\begin{equation}\label{NIRF}
{\mathcal E}(v)=\inf_{j\geq 1}\inf_{\phi\in C(\overline{\Sigma};\Lambda_v^j)}\int_{\Sigma}W\big(\nabla v(x)\mid\phi(x)\big)dx.
\end{equation}
\end{theorem}
\begin{proof}
Fix $v\in \C^1_*(\overline{\Sigma};\RR^3)$ and denote by $\hat{\mathcal E}(v)$ the right-hand side of (\ref{NIRF}). It is easy to verify that ${\mathcal E}(v)\leq\hat{\mathcal E}(v)$. We are thus reduced to prove that
\begin{equation}\label{inequality}
\hat{\mathcal E}(v)\leq{\mathcal E}(v).
\end{equation}
Using Lemma \ref{Lemma3}, we obtain
\begin{equation}\label{inequality1}
\hat{\mathcal E}(v)\leq\inf_{j\geq 1}\int_{\Sigma}\inf_{\zeta\in\Lambda^j_v(x)}W\big(\nabla v(x)\mid\zeta\big)dx.
\end{equation}
Consider the continuous function $\Phi:\overline{\Sigma}\to\RR^3$ defined by 
\begin{equation}\label{ExampleOfContinuousSelection}
\Phi(x):={{\partial_1 v(x)\land\partial_2 v(x)}\over \vert\partial_1v(x)\land\partial_2v(x)\vert^2}.
\end{equation}
 Then, $\det(\nabla v(x)\mid\Phi(x))=1$ for all $x\in\overline{\Sigma}$. Using {\rm (\ref{Icv})} we deduce that there exists $c>0$ depending only on $p$ such that
$$
\int_\Sigma\inf_{\zeta\in\Lambda^{1}_v(x)}W(\nabla v(x)\mid\zeta)dx\leq c\big(|\Sigma|+\|\nabla v\|^p_{L^p(\Sigma;\MM^{3\times 2})}+\|\Phi\|^p_{L^p(\Sigma;\RR^3)}\big).
$$ 
It follows that $\inf_{\zeta\in\Lambda^{1}_v(\cdot)}W(\nabla v(\cdot)\mid\zeta)\in L^1(\Sigma)$. From Lemma \ref{Lemma2}(i) and (ii), we see that $\{\inf_{\zeta\in\Lambda^j_v(\cdot)}W(\nabla v(\cdot)\mid\zeta)\}_{j\geq 1}$ is non-increasing, and that for every  $x\in \overline{\Sigma}$,
\begin{equation}\label{equality}
\inf_{j\geq 1}\inf_{\zeta\in\Lambda^j_v(x)}W\big(\nabla v(x)\mid\zeta\big)=W_0\big(\nabla v(x)\big),
\end{equation}
and (\ref{inequality}) follows from (\ref{inequality1}) and ({\ref{equality}}) by using Lebesgue's dominated convergence theorem.
\end{proof}
           

\section{Existence of $\Gamma\hbox{-}\lim\limits_{\eps\to 0}\mathcal{\I}_\eps$}

In this section we prove Theorem \ref{first_main_result}. Since $\Gamma\hbox{-}\liminf_{\eps\to 0}\mathcal{\I}_\eps\leq\Gamma\hbox{-}\limsup_{\eps\to 0}\mathcal{\I}_\eps$, we only need to show that:
\begin{itemize}
\item[(a)]$\displaystyle\overline{{\mathcal E}}\leq\Gamma\hbox{-}\liminf_{\eps\to 0}\mathcal{\I}_\eps$;
\item[(b)]$\displaystyle\Gamma\hbox{-}\limsup_{\eps\to 0}\mathcal{\I}_\eps\leq\overline{{\mathcal E}}$.
\end{itemize}
In the sequel, we follow the notation used in Section 3.

\subsection{Proof of (a)} 
Let $v\in W^{1,p}(\Sigma;\RR^3)$ and let  $\{v_\eps\}_\eps\subset W^{1,p}(\Sigma;\RR^3)$ be such that $v_\eps\to v$ in $L^{p}(\Sigma;\RR^3)$. We have to prove that
\begin{equation}\label{main_ineq1}
\liminf_{\eps\to 0}{\mathcal E}_\eps(v_\eps)\geq \overline{{\mathcal E}}(v).
\end{equation}
Without loss of generality we can assume that $\sup_{\eps>0}{\mathcal E}_\eps(v_\eps)<+\infty$. To every $\eps>0$ there corresponds $u_\eps\in\pi_{\eps}^{-1}(v_\eps)$ such that 
\begin{equation}\label{main_ineq2}
\E_\eps(v_\eps)\geq E_\eps(u_\eps)-\eps.
\end{equation}
Defining $\hat u_\eps:\Sigma_1\to\RR^3$ by $\hat u_\eps(x,x_3):=u_\eps(x,\eps x_3)$ we have
\begin{equation}\label{main_ineq3}
E_\eps\big(u_\eps\big)=\int_{\Sigma_1}W\Big(\partial_{1} \hat u_\eps(x,x_3)\mid\partial_{2} \hat u_\eps(x,x_3)\mid{1\over\eps}\partial_3 \hat u_\eps(x,x_3)\Big)dxdx_3.
\end{equation}
Using the coercivity of $W$, we deduce that $\left\|{\partial_3 \hat u_\eps}\right\|_{L^p(\Sigma_1;\RR^3)}\le c\eps^{p}$ for all $\eps>0$ and some $c>0$, and so
$
\|\hat u_\eps-v_\eps\|_{L^p(\Sigma_1;\RR^3)}\leq c^\prime\eps^p
$
by Poincar\'e-Wirtinger's inequality, where $c^\prime>0$ is a constant which does not depend on $\eps$. It follows that $\hat u_\eps\to v$ in $L^p(\Sigma_1;\RR^3)$.  For $x_3\in]-{1\over 2},{1\over 2}[$, let $w_\eps^{x_3}\in W^{1,p}(\Sigma;\RR^3)$ given by $w_\eps^{x_3}(x):=\hat u_\eps(x,x_3)$. Then (up to a subsequence) $w_\eps^{x_3}\to v$ in $L^p(\Sigma;\RR^3)$ for a.e. $x_3\in ]-{1\over 2},{1\over 2}[$. Taking (\ref{main_ineq2}) and (\ref{main_ineq3}) into account and using Fatou's lemma, we obtain
$$
\liminf_{\eps\to 0}{{\mathcal E}}_\eps(v_\eps)\geq\int_{-{1\over 2}}^{1\over 2}\left(\liminf_{\eps\to 0}\int_\Sigma W_0\big(\nabla w_\eps^{x_3}(x)\big)dx\right)dx_3,
$$
and so $\liminf_{\eps\to 0}{{\mathcal E}}_\eps(v_\eps)\geq{\mathcal I}(v)$, and (\ref{main_ineq1}) follows by using Theorem \ref{AdDitioNal}.\hfill$\square$

\subsection{Proof of (b)} 
 As $\Gamma\hbox{-}\limsup_{\eps\to 0}\mathcal{\I}_\eps$ is lower semicontinuous with respect to the strong topology of $L^{p}(\Sigma;\RR^3)$ (see \cite[Proposition 6.8 p. 57]{dalmaso0}), it is sufficient to prove that for every $v\in\C^1_*(\overline{\Sigma};\RR^3)$,
\begin{equation}\label{limsupequality}
\limsup_{\eps\to 0}\mathcal{\I}_\eps(v)\leq\mathcal{\I}(v).
\end{equation}  
Given $v\in\C^1_*(\overline{\Sigma};\RR^3)$, fix any $j\geq 1$, and any $n\geq 1$. Using Theorem  \ref{basictheorem} we obtain the existence of $\phi\in C(\overline{\Sigma};\Lambda^j_v)$ such that
\begin{equation}\label{mediainequality}
\int_\Sigma W\big(\nabla v(x)\mid\phi(x)\big)dx\leq{\mathcal E}(v)+{1\over n}.
\end{equation}
By Stone-Weierstrass's approximation theorem, there exists $\{\phi_k\}_{k\geq 1}\subset C^\infty(\overline{\Sigma};\RR^3)$ such that
\begin{equation}\label{uniformityconvergence}
\phi_k\to\phi\hbox{ uniformly as }k\to+\infty.
\end{equation}
We claim that:
\begin{itemize}
\item[(c$_1$)] $\displaystyle\det\big(\nabla v(x)\mid\phi_k(x)\big)\geq{1\over 2j}$ for all $x\in \overline{\Sigma}$, all $k\geq k_v$ and some $k_v\geq 1$;
\item[(c$_2$)] $\displaystyle\lim\limits_{k\to+\infty}\int_\Sigma W\big(\nabla v(x)\mid\phi_k(x)\big)dx=\int_\Sigma W\big(\nabla v(x)\mid\phi(x)\big)dx$.
\end{itemize}
Indeed, setting $\mu_v:=\sup_{x\in\overline{\Sigma}}\vert\partial_1v(x)\land\partial_2v(x)\vert$ ($\mu_v>0$) and using (\ref{uniformityconvergence}), we deduce that there exists $k_v\geq 1$ such that for every $k\geq k_v$,
\begin{equation}\label{supequality}
\sup_{x\in\overline{\Sigma}}\big|\phi_k(x)-\phi(x)\big|<{1\over 2j\mu_v}.
\end{equation}
Let $x\in \overline{\Sigma}$, and let $k\geq k_v$. As $\phi\in C(\overline{\Sigma};\Lambda^j_v)$ we have
\begin{equation}\label{supequality1}
\det\big(\nabla v(x)\mid\phi_k(x)\big)\geq{1\over j}-\det\big(\nabla v(x)\mid\phi_k(x)-\phi(x)\big).
\end{equation}
Noticing that $\det(\nabla v(x)\mid\phi_k(x)-\phi(x))\leq|\partial_1 v(x)\land\partial_2 v(x)||\phi_k(x)-\phi(x)|$, from (\ref{supequality}) and (\ref{supequality1}) we deduce that 
$
\det\big(\nabla v(x)\mid\phi_k(x)\big)\geq{1\over 2j},
$
and (c$_1$) is proved. Combining (c$_1$) with (\ref{Icv}) we see that 
$
\sup_{k\geq k_v}W(\nabla v(\cdot)\mid\phi_k(\cdot))\in L^1(\Sigma).
$ 
As $W$ is continuous we have
$
\lim_{k\to+\infty}W(\nabla v(x)\mid\phi_k(x))=W(\nabla v(x)\mid\phi(x))
$
for all $x\in V$, and (c$_2$) follows by using Lebesgue's dominated convergence theorem, which completes the claim.

Fix any $k\geq k_v$ and define $\theta:]-{1\over 2},{1\over 2}[\to\RR$ by 
$
\theta(x_3):=\inf_{x\in \overline{\Sigma}}\det(\nabla v(x)+x_3\nabla\phi_k(x)\mid\phi_k(x)).
$
Clearly $\theta$ is continuous. By (c$_1$) we have $\theta(0)\geq{1\over 2j}$, and so there exists $\eta_v\in]0,{1\over 2}[$ such that $\theta(x_3)\geq{1\over 4j}$ for all $x_3\in]-\eta_v,\eta_v[$. Let $u_k:\Sigma_1\to\RR$ be given by
$
u_k(x,x_3):=v(x)+x_3\phi_k(x).
$
From the above it follows that
\begin{itemize}
 \item[(c$_3$)] $\det\nabla u_k(x,\eps x_3)\geq{1\over 4j}$ for all $\eps\in]0,\eta_v[$ and all $(x,x_3)\in \overline{\Sigma}\times]-{1\over 2},{1\over 2}[$. 
\end{itemize}
As in the proof of (c$_1$), from (c$_3$) together with (\ref{Icv}) and the continuity of $W$, we obtain
\begin{equation}\label{finalequality}
\lim_{\eps\to 0}{E}_\eps(u_k)=\lim_{\eps\to 0}\int_{\Sigma_1} W\big(\nabla u_k(x,\eps x_3)\big)dxdx_3=\int_\Sigma W\big(\nabla v(x)\mid\phi_k(x)\big)dx.
\end{equation}

For every $\eps>0$ and every $k\geq k_v$, since $\pi_\eps(u_k)=v$ we have ${\mathcal E}_\eps(v)\leq {E}_\eps(u_k)$. Using (\ref{finalequality}), (c$_2$) and (\ref{mediainequality}), we deduce that
$$
\limsup_{\eps\to 0}{\mathcal E}_\eps(v)\leq{\mathcal E}(v)+{1\over n},
$$
and (\ref{limsupequality}) follows by letting $n\to+\infty$.\hfill$\square$


\section{A simplified formula for $\overline{\I}$}

In this section, we prove of Theorem \ref{AdDitioNal}. It is based upon two approximation theorems by Ben Belgacem-Bennequin (see Sect. A.1) and Gromov-Eliasberg (see Sect. A.2). 

Recall the definition of rank one convexity and rank one convex envelope:
\begin{definition} Let $f:\MM^{3\times 2}\to[0,+\infty]$ be a Borel measurable function.
\begin{itemize}
\item[(i)] We say that $f$ is rank one convex if for every $\alpha\in]0,1[$ and every $\xi,\xi^\prime\in\MM^{3\times 2}$ with rank($\xi-\xi^\prime$)=1,
$$
f(\alpha\xi+(1-\alpha)\xi^\prime)\leq \alpha f(\xi)+(1-\alpha)f(\xi^\prime).
$$
\item[(ii)] By the rank one convex envelope of $f$, that we denote by $\mathcal{R}f$, we mean the greatest rank one convex function which less than or equal to $f$.
\end{itemize}
\end{definition}
In \cite[Proposition 7 p. 32 and Lemma 8 p. 34]{benbelgacem-thesis} (see also \cite[Sect. 5.1]{benbelgacem}, \cite[Proposition 3.4.4 p. 112]{trabelsi} and \cite[Lemma 6.5]{trabelsi1}) Ben Belgacem proved the following lemma that we will use in the proof of Theorem \ref{AdDitioNal}. (As $W_0$ is coercive, it is easy to see that ${\mathcal R}W_0$ is coercive.)
\begin{lemma}\label{RW_0Properties}
If {\rm(\ref{Icvbis})} holds then{\rm:}
\begin{itemize}
\item[(i)]  ${\mathcal R}W_0(\xi)\leq c(1+|\xi|^p)$ for all $\xi\in\MM^{3\times 2}$ and some $c>0;$
\item[(ii)] ${\mathcal R}W_0$ is continuous.
\end{itemize}
\end{lemma}
Define $I:W^{1,p}(\Sigma;\RR^3)\to[0,+\infty]$ by
$$
I(v):=\inf\left\{\liminf_{n\to+\infty}\int_\Sigma W_0(\nabla v_n(x))dx:\Aff_{\rm li}(\Sigma;\RR^3)\ni v_n\to v\hbox{ in }L^p(\Sigma;\RR^3)\right\}
$$
with $\Aff_{\rm li}(\Sigma;\RR^3):=\{v\in\Aff(\Sigma;\RR^3):v\hbox{ is locally injective}\}$ ($\Aff(\Sigma;\RR^3)$ is defined in Sect. 2.3). To prove Theorem \ref{AdDitioNal} we will use Proposition \ref{I_li=J_li}.
\begin{proposition}\label{I_li=J_li}
$I=J$ with $J:W^{1,p}(\Sigma;\RR^3)\to[0,+\infty]$ given by
$$
J(v):=\inf\left\{\liminf_{n\to+\infty}\int_\Sigma{\mathcal R} W_0(\nabla v_n(x))dx:\Aff_{\rm li}(\Sigma;\RR^3)\ni v_n\to v\hbox{ in }L^p(\Sigma;\RR^3)\right\}.
$$
\end{proposition}
To prove Proposition \ref{I_li=J_li} we need Lemma \ref{LeMMaBiS} whose proof is contained in the thesis of Ben Belgacem \cite{benbelgacem-thesis}. Since it is difficult to lay hands on this thesis (which is written in French), we give the proof of Lemma \ref{LeMMaBiS} in appendix B.
\begin{lemma}\label{LeMMaBiS} 
$\displaystyle I(v)\leq\int_\Sigma{\mathcal R}W_0(\nabla v(x))dx$ for all $v\in\Aff_{\rm li}(\Sigma;\RR^3)$.
\end{lemma}
\noindent{\em Proof of Proposition {\rm\ref{I_li=J_li}}.} Clearly $J\leq I$. We are thus reduced to prove that 
\begin{equation}\label{inequality}
I\leq J.
\end{equation}
Fix any $v\in W^{1,p}(\Sigma;\RR^3)$ and any sequence $v_n\to v$ in $L^p(\Sigma;\RR^3)$ with $v_n\in\Aff_{\rm li}(\Sigma;\RR^3)$. Using Lemma \ref{LeMMaBiS} we have
$
I(v_n)\leq\int_\Sigma \mathcal{R} W_0(\nabla v_n(x))dx
$
for all $n\geq 1$. Thus,
$$
I(v)\leq\liminf_{n\to+\infty}I(v_n)\leq\liminf_{n\to+\infty}\int_\Sigma \mathcal{R} W_0(\nabla v_n(x))dx,
$$
and (\ref{inequality}) follows.\hfill $\square$

\medskip

\noindent{\em Proof of Theorem {\em\ref{AdDitioNal}}. }We first prove that 
\begin{equation}\label{AddTheorEq1}
\overline{\mathcal{E}}\leq I.
\end{equation}
As in the proof of Proposition \ref{I_li=J_li}, it suffices to show that if $v\in \Aff_{\rm li}(\Sigma;\RR^3)$ then
\begin{equation}\label{AddTheorEq2}
\overline{\mathcal{E}}(v)\leq\int_{\Sigma}W_0(\nabla v(x))dx.
\end{equation}
Let $v\in \Aff_{\rm li}(\Sigma;\RR^3)$. By Theorem \ref{LemmaBBB}-bis (and Lemma \ref{Aff=AffV}), there exists $\{v_n\}_{n\geq 1}\subset C^1_*(\overline{\Sigma};\RR^3)$ such that (\ref{BB_1}) and (\ref{BB_2}) holds and $\nabla v_n(x)\to\nabla v(x)$ a.e. in $\Sigma$. As $W_0$ is continuous (see Lemma \ref{propertiesofW_0}(i)), we have
$$
\lim_{n\to +\infty}W_0\big(\nabla v_n(x)\big)=W_0\big(\nabla v(x)\big)\;\hbox{ a.e. in }\Sigma.
$$
Using (\ref{Icvbis}) together with (\ref{BB_2}), we deduce that there exists $c>0$  such that for every $n\geq 1$ and every measurable set $A\subset\Sigma$,
$$
\int_A W_0\big(\nabla v_n(x)\big)dx\leq c\Big(|A|+\int_A|\nabla v_n(x)-\nabla v(x)|^pdx+\int_A|\nabla v(x)|^pdx\Big).
$$
But $\nabla v_n\to\nabla v$ in $L^p(\Sigma;\MM^{3\times 2})$ by (\ref{BB_1}), hence $\{W_0(\nabla v_n(\cdot))\}_{n\geq 1}$ is absolutely uniformly integrable. Using Vitali's theorem, we obtain
$$
\lim_{n\to+\infty}\int_{\Sigma}W_0(\nabla v_n(x))dx=\int_{\Sigma}W_0(\nabla v(x))dx,
$$ 
and (\ref{AddTheorEq2}) follows.

We now prove that
\begin{equation}\label{AddTheorEq3}
J\leq\overline{J},
\end{equation}
with $\overline{J}:W^{1,p}(\Sigma;\RR^3)\to[0,+\infty]$ given by
$$
\overline{J}(v):=\inf\left\{\liminf_{n\to+\infty}\int_\Sigma{\mathcal R} W_0(\nabla v_n(x))dx:W^{1,p}(\Sigma;\RR^3)\ni v_n\to v\hbox{ in }L^p(\Sigma;\RR^3)\right\}.
$$
It is sufficient to show that
\begin{equation}\label{AddTheorEq4}
J(v)\leq\int_{\Sigma}\mathcal{R}W_0(\nabla v(x))dx.
\end{equation} 
Let $v\in W^{1,p}(\Sigma;\RR^3)$. By Corollary \ref{GromovEliasbergConsequence}, there exists $\{v_n\}_{n\geq 1}\subset\Aff_{\rm li}(\Sigma;\RR^3)$ such that $\nabla v_n\to\nabla v$ in $L^p(\Sigma;\RR^3)$ and  $\nabla v_n(x)\to\nabla v(x)$ a.e. in $\Sigma$. Taking Lemma \ref{RW_0Properties} into account, from Vitali's lemma, we see that
$$
\lim_{n\to+\infty}\int_{\Sigma}\mathcal{R}W_0(\nabla v_n(x))dx=\int_{\Sigma}\mathcal{R}W_0(\nabla v(x))dx,
$$
and (\ref{AddTheorEq4}) follows.

Noticing that $\mathcal{I}\leq\overline{\mathcal{E}}$ and $\overline{J}\leq\mathcal{I}$, and combining Proposition \ref{I_li=J_li} with (\ref{AddTheorEq1}) and (\ref{AddTheorEq3}), we conclude that $\overline{\mathcal{E}}=\mathcal{I}$.\hfill$\square$


\appendix

\section{Approximation theorems}

\subsection{Ben Belgacem-Bennequin's theorem} Denote by $\Aff^{ET}(\Sigma;\RR^3)$ the space of Ekeland-Temam continuous piecewise affine functions from $\Sigma$ to $\RR^3$, i.e., {\em $u\in\Aff^{ET}(\Sigma;\RR^3)$ if and only if $v$ is continuous and there exists a finite family $(V_i)_{i\in I}$ of open disjoint subsets of $\Sigma$ such that $|\Sigma\setminus \cup_{i\in I} V_i|=0$ and for every $i\in I$, the restriction of $v$ to $V_i$ is affine.} Note that from Ekeland-Temam \cite{ekeland}, we know that $\Aff^{ET}(\Sigma;\RR^3)$ is strongly dense in $W^{1,p}(\Sigma;\RR^3)$. Set
$$
\Aff^{ET}_{\rm li}(\Sigma;\RR^3):=\Big\{v\in\Aff^{ET}(\Sigma;\RR^3):v\hbox{ is locally injective}\Big\}.
$$
In  \cite[Lemma 8 p. 114]{benbelgacem-thesis} (see also \cite[Proposition C.0.4 p. 127]{trabelsi} and \cite[Lemma 1.3]{trabelsi1}) Ben Belgacem and Bennequin proved the following result. 
\begin{theorem}\label{LemmaBBB}
For every $v\in\Aff^{ET}_{\rm li}(\Sigma;\RR^3)$, there exists $\{v_n\}_{n\geq 1}\subset C^1_*(\overline{\Sigma};\RR^3)$ such that{\rm:} 
\begin{equation}\label{BB_1}
\hbox{$v_n\to v$ {\em in } $W^{1,p}(\Sigma;\RR^3)${\rm;}}
\end{equation}
\begin{equation}\label{BB_2}
\hbox{$|\partial_1v_n(x)\land\partial_2v_n(x)|\geq\delta$ for all $x\in\overline{\Sigma}$, all $n\geq 1$ and some $\delta>0$.} 
\end{equation}
\end{theorem}
Denote by $\Aff^V(\Sigma;\RR^3)$ the space of Vitali continuous piecewise affine functions from $\Sigma$ to $\RR^3$ (introduced by Ben Belgacem in \cite{benbelgacem-thesis,benbelgacem}), i.e., {\em $v\in\Aff^V(\Sigma;\RR^3)$ if and only if $v$ is continuous and there exists a finite or countable family $(O_i)_{i\in I}$ of dsjoint open subsets of $\Sigma$ such that $|\partial O_i|=0$ for all $i\in I$, $|\Sigma\setminus\cup_{i\in I}O_i|=0$, and $v(x)=\xi_i\cdot x+a_i$ if $x\in O_i$, where $a_i\in\RR^3$, $\xi_i\in\MM^{3\times 2}$ and ${\rm Card}\{\xi_i:i\in I\}$ is finite.} In \cite[Lemma 3.1.5 p. 99]{trabelsi} Trabelsi remarked that Theorem \ref{LemmaBBB} can be generalized replacing the space $\Aff^{ET}_{\rm li}(\Sigma;\RR^{3})$ by
$$
\Aff^V_{\rm li}(\Sigma;\RR^3):=\Big\{v\in\Aff^V(\Sigma;\RR^3):v\hbox{ is locally injective}\Big\}.
$$

\noindent{\bf Theorem \ref{LemmaBBB}-bis.} {\em For every $v\in\Aff^{V}_{\rm li}(\Sigma;\RR^3)$, there exists $\{v_n\}_{n\geq 1}\subset C^1_*(\overline{\Sigma};\RR^3)$ satisfying} (\ref{BB_1}) and (\ref{BB_2}).

\medskip

Here we consider the space $\Aff(\Sigma;\RR^3)$ defined in Sect. 2.3. It is clear that $\Aff^{ET}(\Sigma;\RR^3)\subset\Aff(\Sigma;\RR^3)$, and so $\Aff(\Sigma;\RR^3)$ is strongly dense in $W^{1,p}(\Sigma;\RR^3)$. Moreover, we have
\begin{lemma}\label{Aff=AffV}
$\Aff^V(\Sigma;\RR^3)=\Aff(\Sigma;\RR^3)$. 
\end{lemma}
\begin{proof}
Setting $D_i:=\{x\in\cup_{i\in I}O_i:\nabla v(x)=\xi_i\}$ with $v\in\Aff^V(\Sigma;\RR^3)$,  we see that ${\rm Card}\{D_i:i\in I\}$ is finite. Thus $\Aff^V(\Sigma;\RR^3)\subset\Aff(\Sigma;\RR^3)$. Given $v\in\Aff(\Sigma;\RR^3)$, let $(O_j)_{j\in J_i}$ be the connected components of $D_i$ with $i\in I$ (where $I$ is finite). Since $D_i$ is open, $O_j$ is open for all $j\in J_i$, hence $J_i$ is finite or countable because $\QQ^2$ is dense in $\RR^2$. Moreover, for each $j\in J_i$, the restriction of $v$ to $O_j$ is affine. Thus $\Aff(\Sigma;\RR^3)\subset\Aff^V(\Sigma;\RR^3)$.
\end{proof}

\subsection{Gromov-Eliashberg's theorem} In \cite[Theorem 1.3.4B]{gromovEliashberg} (see also \cite[Theorem B$^\prime_1$ p. 20]{gromov}) Gromov and Eliashberg proved the following result. 
\begin{theorem}
Let $1\leq N<m$ be two integers and let $M$ be a compact $N$-di-mensional manifold which can be immersed in $\RR^m$. Then, for each $C^1$-differentiable function $v$ from $M$ to $\RR^m$ there exists a sequence $\{v_n\}_n$ of $C^1$-immersions from $M$ to $\RR^m$ such that $v_n\to v$ in $W^{1,p}(M;\RR^m)$.
\end{theorem}
In our context, we have
\begin{theorem}\label{GromovEliasberg}
For every $v\in C^1(\overline{\Sigma};\RR^3)$ there exists $\{v_n\}_{n\geq 1}\subset C^1_*(\overline{\Sigma};\RR^3)$ such that $v_n\to v$ in $W^{1,p}(\Sigma;\RR^3)$.
\end{theorem}
Moreover, from \cite[Proposition 3.1.7 p. 100]{trabelsi}, we have
\begin{proposition}\label{TrabelsiThesis}
For every $v\in C^1_*(\overline{\Sigma};\RR^3)$ there exists $\{v_n\}_{n\geq 1}\subset\Aff^{ET}_{\rm li}(\Sigma;\RR^3)$ such that  $v_n\to v$ in $W^{1,p}(\Sigma;\RR^3)$.
\end{proposition}
Thus, as a consequence of Theorem \ref{GromovEliasberg} and Proposition \ref{TrabelsiThesis}, we obtain
\begin{corollary}\label{GromovEliasbergConsequence}
$\Aff^{ET}_{\rm li}(\Sigma;\RR^3)$ is strongly dense in $W^{1,p}(\Sigma;\RR^3)$.
\end{corollary}


\section{Ben Belgacem's lemma}

In this appendix we prove Ben Belgacem's lemma, i.e., Lemma \ref{LeMMaBiS}.

\subsection{Preliminaries.} Define the sequence $\{\R_i W_0\}_{i\geq 0}$ by $\R_0 W_0=W_0$ and for every $i\geq 1$ and every $\xi\in\MM^{3\times 2}$,
$$
\R_{i+1}W_0(\xi):=\infff\limits_{t\in[0,1]}\Big\{(1-t)\R_i W_0(\xi-t a\otimes b)+t\R_i W_0(\xi+(1-t)a\otimes b)\Big\}.
$$
Recall that $W_0$ is coercive and continuous (see Lemma \ref{propertiesofW_0}(i)). The following lemma is due to Kohn and Strang \cite{kohnstrang}.
\begin{lemma}\label{Khon-Strang}
$\R_{i+1} W_0\leq \R_iW_0$ for all $i\geq 0$ and $\R W_0=\inf_{i\geq 0} \R_i W_0$.
\end{lemma}
Fix any $i\geq 0$ and any $v\in\Aff_{\rm li}(\Sigma;\RR^3):=\{v\in\Aff(\Sigma;\RR^3):v\hbox{ is locally injective}\}$ (with $\Aff(\Sigma;\RR^3)$ defined in Sect. 2.3). By definition, there exists a finite family $(V_j)_{j\in J}$ of open disjoint subsets of $\Sigma$ such that $|\partial V_j|=0$ for all $j\in J$, $|\Sigma\setminus\cup_{j\in J}V_j|=0$ and, for every $j\in J$, $\nabla v(x)=\xi_j$ in $V_j$ with $\xi_j\in\MM^{3\times 2}$. (As $v$ is locally injective we have ${\rm rank}(\xi_j)=2$ for all $j\in J$.) Fix any $j\in J$. For a proof of Lemmas \ref{BBLemma2} and \ref{BBLemma1} we refer to \cite[Proposition 3.1.2 p. 96]{trabelsi}. 
\begin{lemma}\label{BBLemma2}
$\R_i W_0$ is continuous.
\end{lemma}
\begin{lemma}\label{BBLemma1} 
There exist $a\in\RR^2$, $b\in\RR^3$ and $t\in[0,1]$ such that
$$
\R_{i+1}W_0(\xi_j)=(1-t)\R_i W_0(\xi_j-t a\otimes b)+t\R_i W_0(\xi_j+(1-t)a\otimes b).
$$
\end{lemma}
Without loss of generality we can assume that $a=(1,0)$. For each $n\geq 3$ and  each $k\in\{0,\cdots,n-1\}$, consider $A^-_{k,n},A^+_{k,n},B_{k,n},B^-_{k,n},B^+_{k,n},C_{k,n},C^-_{k,n},C^+_{k,n}\subset Y$ given by:

\smallskip

$A^-_{k,n}:=\big\{(x_1,x_2)\in Y:{k\over n}\leq x_1\leq{k\over n}+{1-t\over n}\hbox{ and }{1\over n}\leq x_2\leq 1-{1\over n}\big\}$;

$A^+_{k,n}:=\big\{(x_1,x_2)\in Y:{k\over n}+{1-t\over n}\leq x_1\leq {k+1\over n}\hbox{ and }{1\over n}\leq x_2\leq 1-{1\over n}\big\}$;

$B_{k,n}:=\big\{(x_1,x_2)\in Y:{k\over n}\leq x_1\leq{k+1\over n}\hbox{ and }0\leq x_2\leq -x_1+{k+1\over n}\big\}$;

$B^-_{k,n}:=\big\{(x_1,x_2)\in Y:-x_2+{k+1\over n}\leq x_1\leq-tx_2+{k+1\over n}\hbox{ and }0\leq x_2\leq {1\over n}\big\}$;

$B^+_{k,n}:=\big\{(x_1,x_2)\in Y:-tx_2+{k+1\over n}\leq x_1\leq {k+1\over n}\hbox{ and }0\leq x_2\leq {1\over n}\big\}$;

$C_{k,n}:=\big\{(x_1,x_2)\in Y:{k\over n}\leq x_1\leq{k+1\over n}\hbox{ and }x_1+1-{k+1\over n}\leq x_2\leq 1\big\}$;

$C^-_{k,n}:=\big\{(x_1,x_2)\in Y:x_2-1+{k+1\over n}\leq x_1\leq t(x_2-1)+{k+1\over n}\hbox{ and }0\leq x_2\leq {1\over n}\big\}$;

$C^+_{k,n}:=\big\{(x_1,x_2)\in Y:t(x_2-1)+{k+1\over n}\leq x_1\leq {k+1\over n}\hbox{ and }0\leq x_2\leq {1\over n}\big\}$,

\smallskip

\noindent  and define $\{\sigma_{n}\}_{n\geq 1}\subset\Aff_0(Y;\RR)$ by
$$
\sigma_{n}(x_1,x_2):=\left\{
\begin{array}{ll}
-t(x_1-{k\over n})&\hbox{if }(x_1,x_2)\in A^-_{k,n}\\
(1-t)(x_1-{k+1\over n})&\hbox{if }(x_1,x_2)\in A^+_{k,n}\cup B^+_{k,n}\cup C^+_{k,n}\\
-t(x_1+x_2-{k+1\over n})&\hbox{if }(x_1,x_2)\in B^-_{k,n}\\
-t(x_1-x_2+1-{k+1\over n})&\hbox{if }(x_1,x_2)\in C^-_{k,n}\\
0&\hbox{if }(x_1,x_2)\in B_{k,n}\cup C_{k,n}
\end{array}
\right.
$$
(see Figure B.1).

\smallskip

\begin{picture}(300,255)
\put(-5,40){\vector(1,0){210}}
\put(195,34){\SMALL$x_1$}
\put(-10,240){\SMALL$x_2$}
\put(0,39.9){\line(1,0){180}}
\put(0,40.1){\line(1,0){180}}
\put(0,40){\circle*{3}}
\put(-5,34){\SMALL$0$}
\put(20,40){\circle*{3}}
\put(12,31){\SMALL$1-t\over n$}
\put(30,40){\circle*{3}}
\put(26.5,31){\SMALL$1\over n$}
\put(50,40){\circle*{3}}
\put(42,31){\SMALL${2-t\over n}$}
\put(60,40){\circle*{3}}
\put(56.5,31){\SMALL$2\over n$}
\put(71,32){\SMALL$\cdots$}
\put(90,40){\circle*{3}}
\put(86.5,31){\SMALL$k\over n$}
\put(120,40){\circle*{3}}
\put(115,31){\SMALL$k+1\over n$}
\put(93,31){\Tiny$k+1-t\over n$}
\put(110,40){\circle*{3}}
\put(131.5,32){\SMALL$\cdots$}
\put(150,40){\circle*{3}}
\put(144,31){\SMALL$n-1\over n$}
\put(170,40){\circle*{3}}
\put(180,40){\circle*{3}}
\put(162,31){\SMALL$n-t\over n$}
\put(178,32){\SMALL$1$}
\put(0,35){\vector(0,1){210}}
\put(0,220){\circle*{3}}
\put(-5,218){\SMALL$1$}
\put(0,190){\circle*{3}}
\put(-17,188){\SMALL$n-1\over n$}
\put(0,70){\circle*{3}}
\put(-8,68){\SMALL$1\over n$}
\put(0,220){\line(1,0){180}}
\put(0,220.1){\line(1,0){180}}
\put(0,219.9){\line(1,0){180}}
\put(0,190){\line(1,0){20}}
\put(30,190){\line(1,0){20}}
\put(60,190){\line(1,0){20}}
\put(90,190){\line(1,0){20}}
\put(120,190){\line(1,0){20}}
\put(150,190){\line(1,0){20}}
\put(150,70){\line(1,0){20}}
\put(120,70){\line(1,0){20}}
\put(90,70){\line(1,0){20}}
\put(60,70){\line(1,0){20}}
\put(30,70){\line(1,0){20}}
\put(0,70){\line(1,0){20}}
\put(-0.1,40){\line(0,1){180}}
\put(0.1,40){\line(0,1){180}}
\put(180,40){\line(0,1){180}}
\put(180.1,40){\line(0,1){180}}
\put(179.9,40){\line(0,1){180}}
\put(150,40){\line(0,1){180}}
\put(150.1,40){\line(0,1){180}}
\put(149.9,40){\line(0,1){180}}
\put(120,40){\line(0,1){180}}
\put(120.1,40){\line(0,1){180}}
\put(119.9,40){\line(0,1){180}}
\put(90,40){\line(0,1){180}}
\put(90.1,40){\line(0,1){180}}
\put(89.9,40){\line(0,1){180}}
\put(60,40){\line(0,1){180}}
\put(60.1,40){\line(0,1){180}}
\put(59.9,40){\line(0,1){180}}
\put(30,40){\line(0,1){180}}
\put(30.1,40){\line(0,1){180}}
\put(29.9,40){\line(0,1){180}}
\put(170,70){\line(0,1){120}}
\put(140,70){\line(0,1){120}}
\put(110,70){\line(0,1){120}}
\put(80,70){\line(0,1){120}}
\put(50,70){\line(0,1){120}}
\put(20,70){\line(0,1){120}}
\put(0,70){\line(1,-1){30}}
\put(30,70){\line(1,-1){30}}
\put(60,70){\line(1,-1){30}}
\put(90,70){\line(1,-1){30}}
\put(120,70){\line(1,-1){30}}
\put(150,70){\line(1,-1){30}}
\put(0,190){\line(1,1){30}}
\put(30,190){\line(1,1){30}}
\put(60,190){\line(1,1){30}}
\put(90,190){\line(1,1){30}}
\put(120,190){\line(1,1){30}}
\put(150,190){\line(1,1){30}}
\put(20,70){\line(1,-3){10}}
\put(50,70){\line(1,-3){10}}
\put(80,70){\line(1,-3){10}}
\put(110,70){\line(1,-3){10}}
\put(140,70){\line(1,-3){10}}
\put(170,70){\line(1,-3){10}}
\put(20,190){\line(1,3){10}}
\put(50,190){\line(1,3){10}}
\put(80,190){\line(1,3){10}}
\put(110,190){\line(1,3){10}}
\put(140,190){\line(1,3){10}}
\put(170,190){\line(1,3){10}}
\put(260,40){\line(1,0){30}}
\put(260,40.1){\line(1,0){30}}
\put(260,39.9){\line(1,0){30}}
\put(235,52){\Tiny$B_{k,n}$}\put(252,52){\vector(1,0){17}}
\put(235,208){\Tiny$C_{k,n}$}\put(252,208){\vector(1,0){17}}
\put(235,128){\Tiny$A^-_{k,n}$}\put(252,128){\vector(1,0){20}}
\put(307,128){\Tiny$A^+_{k,n}$}\put(304,128){\vector(-1,0){20}}
\put(307.5,198){\Tiny$C^+_{k,n}$}\put(304.5,198){\vector(-1,0){20}}
\put(307.5,62){\Tiny$B^+_{k,n}$}\put(304.5,62){\vector(-1,0){20}}
\put(235,62){\Tiny$B^-_{k,n}$}\put(252,62){\vector(1,0){25}}
\put(235,198){\Tiny$C^-_{k,n}$}\put(252,198){\vector(1,0){25}}
\put(260,220){\line(1,0){30}}
\put(260,220.1){\line(1,0){30}}
\put(260,219.9){\line(1,0){30}}
\put(260,40){\circle*{3}}
\put(256.5,31){\SMALL$k\over n$}
\put(290,40){\circle*{3}}
\put(285,31){\SMALL$k+1\over n$}
\put(263,31){\Tiny$k+1-t\over n$}
\put(280,40){\circle*{3}}
\put(290,40){\line(0,1){180}}
\put(290.1,40){\line(0,1){180}}
\put(289.9,40){\line(0,1){180}}
\put(260,40){\line(0,1){180}}
\put(260.1,40){\line(0,1){180}}
\put(259.9,40){\line(0,1){180}}
\put(260,70){\line(1,-1){30}}
\put(280,70){\line(1,-3){10}}
\put(280,190){\line(1,3){10}}
\put(260,190){\line(1,1){30}}
\put(280,70){\line(0,1){120}}
\put(260,70){\line(1,0){30}}
\put(260,190){\line(1,0){30}}
\put(0,10){\SMALL Figure B.1. The function $\sigma_n$ and the sets $A^-_{k,n},A^+_{k,n},B_{k,n},B^-_{k,n},B^+_{k,n},C_{k,n},C^-_{k,n},C^+_{k,n}$.}
\end{picture}

\noindent Set
$$
b_{\ell}:=\left\{
\begin{array}{ll}
b&\hbox{if }b\not\in{\rm Im}\xi_j\\
b+{1\over\ell}\nu&\hbox{if }b\in{\rm Im}\xi_j
\end{array}
\right.
$$
(with ${\rm Im}\xi_j:=\{\xi_j\cdot x:x\in\RR^2\}$) where $\ell\geq 1$ and $\nu\in\RR^3$ is a normal vector to ${\rm Im}\xi_j$.
\begin{lemma}\label{LeMMaBBFund2}
Define $\{\theta_{n,\ell}\}_{n,\ell\geq 1}\subset\Aff_0(Y;\RR^3)$ by 
$$
\theta_{n,\ell}(x):=\sigma_{n}(x)b_{\ell}.
$$
Then{\rm:}
\begin{itemize}
\item[(i)] for every $\ell\geq 1$, $\theta_{n,\ell}\to 0$ in $L^p(Y;\RR^3);$
\item[(ii)] $\displaystyle\lim_{\ell\to+\infty}\lim_{n\to+\infty}\int_Y\R_iW_0(\xi_j+\nabla \theta_{n,\ell}(x))dx=\R_{i+1}W_0(\xi_j)$.
\end{itemize}
\end{lemma}
\begin{proof}
(i) It suffices to prove that $\sigma_n\to 0$ in $L^p(Y;\RR)$.  For every $k\in\{0,\cdots,n-1\}$, it is clear that $|\sigma_n(x)|^p\leq {t^p(1-t)^p\over n^p}$ for all $x\in]{k\over n},{k+1\over n}[\times]0,1[$, and so
$$
\int_{]{k\over n},{k+1\over n}[\times]0,1[}|\sigma_n(x)|^pdx\leq {t^p(1-t)^p\over n^{p+1}}.
$$
As
$$
\int_Y|\sigma_{n}(x)|^pdx=\sum_{k=0}^{n-1}\int_{]{k\over n},{k+1\over n}[\times]0,1[}|\sigma_{n}(x)|^pdx
$$
it follows that 
$$
\int_Y|\sigma_{n}(x)|^pdx\leq{t^p(1-t)^p\over n^{p}},
$$
which gives the desired conclusion.

(ii) Recalling that $a=(1,0)$ we see that
$$
\xi_j+\nabla\theta_{n,\ell}(x):=\left\{
\begin{array}{ll}
\xi_j-ta\otimes b_\ell&\hbox{if }x\in {\rm int}(A^-_{k,n})\\
\xi_j+(1-t)a\otimes b_\ell&\hbox{if }x\in {\rm int}(A^+_{k,n}\cup B^+_{k,n}\cup C^+_{k,n})\\
\xi_j-t(a+a^\perp)\otimes b_\ell&\hbox{if }x\in {\rm int}(B^-_{k,n})\\
\xi_j-t(a-a^\perp)\otimes b_\ell&\hbox{if }x\in {\rm int}(C^-_{k,n})\\
\xi_j&\hbox{if }x\in {\rm int}(B_{k,n})\cup{\rm int}(C_{k,n})
\end{array}
\right.
$$
with $a^\perp=(0,1)$ (and ${\rm int}(E)$ denotes the interior of the set $E$). Moreover, we have:

\medskip

${\displaystyle\int_{\cup_{k=0}^{n-1} A^-_{k,n}}\R_iW_0(\xi_j-ta\otimes b_\ell)dx}=(1-t)(1-{2\over n})\R_iW_0(\xi_j-ta\otimes b_\ell)$;

${\displaystyle\int_{\cup_{k=0}^{n-1} A^+_{k,n}}\R_i W_0(\xi_j+(1-t)a\otimes b_\ell)dx}=t(1-{2\over n})\R_i W_0(\xi_j+(1-t)a\otimes b_\ell)$;

${\displaystyle\int_{\cup_{k=0}^{n-1} (B^+_{k,n}\cup C^+_{k,n})}\R_i W_0(\xi_j+(1-t)a\otimes b_\ell)dx}={t\over n}\R_i W_0(\xi_j+(1-t)a\otimes b_\ell)$;

${\displaystyle\int_{\cup_{k=0}^{n-1} B^-_{k,n}}\R_i W_0(\xi_j-t(a+a^\perp)\otimes b_\ell)dx}={1-t\over 2n}\R_i W_0(\xi_j-t(a+a^\perp)\otimes b_\ell)$;

${\displaystyle\int_{\cup_{k=0}^{n-1} C^-_{k,n}}\R_i W_0(\xi_j-t(a-a^\perp)\otimes b_\ell)dx}={1-t\over 2n}\R_i W_0(\xi_j-t(a-a^\perp)\otimes b_\ell)$;

${\displaystyle\int_{\cup_{k=0}^{n-1} (B_{k,n}\cup C_{k,n})}\R_i W_0(\xi_j)dx}={1\over n}\R_i W_0(\xi_j)$.

\medskip

\noindent Hence
\begin{eqnarray*}
\int_Y\R_i W_0(\xi_j+\nabla\theta_{n,\ell}(x))dx&=&\Big(1-{2\over n}\Big)\Big[(1-t)\R_iW_0(\xi_j-ta\otimes b_\ell)+t\R_i W_0(\xi_j\\
&&+(1-t)a\otimes b_\ell)\Big]+{1\over n}\Big[t\R_i W_0(\xi_j+(1-t)a\otimes b_\ell)\\
&&+{1-t\over 2}\big(\R_i W_0(\xi_j-t(a+a^\perp)\otimes b_\ell)+\R_i W_0(\xi_j\hskip-0.8mm-\\
&&t(a-a^\perp)\otimes b_\ell)\big)+\R_i W_0(\xi_j)\Big]
\end{eqnarray*}
for all $n,\ell\geq 1$. It follows that for every $\ell\geq 1$,
\begin{eqnarray*}
\lim_{n\to+\infty}\int_Y\R_i W_0(\xi_j+\nabla\theta_{n,\ell}(x))dx&=&(1-t)\R_iW_0(\xi_j-ta\otimes b_\ell)\\
&&+t\R_i W_0(\xi_j+(1-t)a\otimes b_\ell).
\end{eqnarray*}
Taking Lemma \ref{BBLemma2} into account and noticing that $b_\ell\to b$, we deduce that
\begin{eqnarray*}
\lim_{\ell\to+\infty}\lim_{n\to+\infty}\int_Y\R_i W_0(\xi_j+\nabla\theta_{n,\ell}(x))dx&=&(1-t)\R_iW_0(\xi_j-ta\otimes b)\\
&&+t\R_i W_0(\xi_j+(1-t)a\otimes b),
\end{eqnarray*}
and (ii) follows by using Lemma \ref{BBLemma1}.
 \end{proof}
Consider $V^j_{q}\subset V_j$ given by $V^j_{q}:=\{x\in V_j:{\rm dist}(x,\partial V_j)>{1\over q}\}$ with $q\geq 1$ large enough. By Vitali's covering theorem, there exists a finite or countable family $(r_{m}+\rho_{m}Y)_{m\in M}$ of disjoint subsets of $V^j_{q}$, with $r_{m}\in\RR^2$ and $\rho_{m}\in]0,1[$, such that
$
|V^j_{q}\setminus\cup_{m\in M}(r_{m}+\rho_{m}Y)|=0
$
 (and so $\sum_{m\in M}\rho_{m}^2=|V^j_{q}|$). Let $\{\phi_{n,\ell,q}\}_{n,\ell,q\geq 1}\subset\Aff_0(V_j;\RR^3)$ be given by
$$
\phi_{n,\ell,q}(x):=\left\{
\begin{array}{ll}
\displaystyle\rho_{m}\theta_{n,\ell}\left({x-r_{m}\over \rho_{m}}\right)&\hbox{ if }x\in r_{m}+\rho_{m}Y\subset V^j_{q}\\
0&\hbox{if }x\in V_j\setminus V^j_{q}.
\end{array}
\right.
$$
 \begin{lemma}\label{LeMMaBBFund3}
 Define $\{\Phi^j_{n,\ell,q}\}_{n,\ell,q\geq 1}\subset\Aff(V_j;\RR^3)$ by
\begin{equation}\label{BBFunct2}
\Phi^j_{n,\ell,q}(x):=v(x)+\phi_{n,\ell,q}(x).
\end{equation}
Then{\rm:}
\begin{itemize}
\item[(i)] for every $n,\ell,q\geq 1$, $\Phi^j_{n,\ell,q}$ is locally injective{\rm;}
\item[(ii)] for every $\ell,q\geq 1$, $\Phi^j_{n,\ell,q}\to v$ in $L^p(V_j;\RR^3);$
\item[(iii)] $\displaystyle\lim_{q\to+\infty}\lim_{\ell\to+\infty}\lim_{n\to+\infty}\int_{V_j}\R_iW_0(\nabla \Phi^j_{n,\ell,q}(x))dx=|V_j|\R_{i+1}W_0(\xi_j)$.
\end{itemize}
\end{lemma}
\begin{proof}
(i) Let $x\in V_j$ and let $W\subset V_j$ be the connected component of $V_j$ such that $x\in W$ (as $V_j$ is open, so is $W$). Since $\nabla v=\xi_j$ in $W$, there exists $c\in\RR^3$ such that $v(x^\prime)=\xi_j\cdot x^\prime+c$ for all $x^\prime\in W$.   We claim that ${\Phi^j_{n,\ell,q}}{\lfloor_{W}}$ is injective. Indeed, let $x^\prime\in W$ be such that $\Phi^j_{n,\ell,q}(x)=\Phi^j_{n,\ell,q}(x^\prime)$. One the three possibilities holds: 
\begin{itemize}
\item[(a)] $\Phi^j_{n,\ell,q}(x)=\xi_j\cdot x+c+\rho_m\sigma_{n}\big({x-r_m\over\rho_m}\big)b_\ell$ and $\Phi^j_{n,\ell,q}(x^\prime)=\xi_j\cdot x^\prime+c+\rho_{m^\prime}\sigma_{n}\big({x^\prime-r_{m^\prime}\over\rho_{m^\prime}}\big)b_\ell$;  
\item[(b)] $\Phi^j_{n,\ell,q}(x)=\xi_j\cdot x+c+\rho_m\sigma_{n,\ell}\big({x-r_m\over\rho_m}\big)b_\ell$ and $\Phi^j_{n,\ell,q}(x^\prime)=\xi_j\cdot x^\prime+c$;
\item[(c)] $\Phi^j_{n,\ell,q}(x)=\xi_j\cdot x+c$ and $\Phi^j_{n,\ell,q}(x^\prime)=\xi_j\cdot x^\prime+c$.
\end{itemize}
Setting $\alpha:=\rho_m\sigma_n({x-r_m\over\rho_m})-\rho_{m^\prime}\sigma_n({x^\prime-r_{m^\prime}\over \rho_{m^\prime}})$ and $\beta:=\rho_m\sigma_n({x-r_m\over\rho_m})$ we have:
$$
\left\{
\begin{array}{ll}
\xi_j(x^\prime-x)=0&\hbox{if }\alpha=0\\
b_\ell={1\over \alpha}\xi_j(x^\prime-x)&\hbox{if }\alpha\not=0
\end{array}
\right.
\hbox{ when (a) is satisfied;}
$$
$$
\left\{
\begin{array}{ll}
\xi_j(x^\prime-x)=0&\hbox{if }\beta=0\\
b_\ell={1\over \beta}\xi_j(x^\prime-x)&\hbox{if }\beta\not=0
\end{array}
\right.
\hbox{ when (b) is satisfied;}
$$
$$
\xi_j(x^\prime-x)=0\hbox{ when (c) is satisfied.}
$$
It follows that if $x\not=x^\prime$ then either ${\rm rank}(\xi_j)<2$ or $b_\ell\in{\rm Im}\xi_j$ which is impossible. Hence $x=x^\prime$, and the claim is proved. Thus $\Phi^j_{n,\ell,q}$ is locally injective.

(ii) As $\rho_m\in]0,1[$ for all $m\in M$ and $\sum_{m\in M}\rho_m^2=|V^j_q|$ we have
$$
\int_{V^j_q}|\phi_{n,\ell,q}(x)|^pdx\leq|V^j_q|\int_Y|\theta_{n,\ell}(x)|^pdx.
$$
Using Lemma \ref{LeMMaBBFund2}(i) we deduce that for every $\ell,q\geq 1$,
$$
\lim_{n\to+\infty}\int_{V^j_q}|\phi_{n,\ell,q}(x)|^pdx=0,
$$
and (ii) follows.

(iii) Recalling that $\sum_{m\in M}\rho_m^2=|V^j_q|$ we see that
\begin{eqnarray*}
\int_{V_j}\R_iW_0(\nabla \Phi^j_{n,\ell,q}(x))dx\hskip-2mm&=&\hskip-2mm\int_{V_j}\R_iW_0(\xi_j+\nabla \phi_{n,\ell,q}(x))dx\\
&=&\hskip-2mm\int_{V^j_q}\R_iW_0(\xi_j+\nabla \phi_{n,\ell,q}(x))dx+|V_j\setminus V^j_q|\R_iW_0(\xi_j)\\
&=&\hskip-2mm|V^j_q|\int_Y\R_iW_0(\xi_j+\nabla\theta_{n,\ell}(x))dx+|V_j\setminus V^j_q|\R_iW_0(\xi_j).
\end{eqnarray*}
Using Lemma \ref{LeMMaBBFund2}(ii) we deduce that for every $q\geq 1$,
$$
\lim_{\ell\to+\infty}\lim_{n\to+\infty}\int_{V_j}\R_iW_0(\nabla \Phi^j_{n,\ell,q}(x))dx=|V^j_q|\R_{i+1}W_0(\xi_j)+|V_j\setminus V^j_q|\R_iW_0(\xi_j),
$$
and (iii) follows by noticing that $|V^j_q|\to|V_j|$ and $|V_j\setminus V^j_q|\to 0$.
\end{proof}

\subsection{Proof of Lemma \ref{LeMMaBiS}} According to Lemma \ref{Khon-Strang}, it is sufficient to show that for every $i\geq 0$,
$$
I(v)\leq \int_\Sigma\R_iW_0(\nabla v(x))dx\hbox{ for all }v\in\Aff_{\rm li}(\Sigma;\RR^3).\leqno (P_i)
$$
The proof is by induction on $i$. As $R_0 W_0=W_0$ it is clear that $(P_0)$ is true. Assume that $(P_i)$ is true, and prove that $(P_{i+1})$ is true. Let $v\in\Aff_{\rm li}(\Sigma;\RR^3)$. By definition, there exists a finite family $(V_j)_{j\in J}$ of open disjoint subsets of $\Sigma$ such that $|\partial V_j|=0$ for all $j\in J$, $|\Sigma\setminus\cup_{j\in J}V_j|=0$ and, for every $j\in J$, $\nabla v(x)=\xi_j$ in $V_j$ with $\xi_j\in\MM^{3\times 2}$. Define $\{\Psi_{n,\ell,q}\}_{n,\ell,q\geq 1}\subset\Aff(\Sigma;\RR^3)$ by
$$
\Psi_{n,\ell,q}(x):=\Phi^j_{n,\ell,q}(x)\hbox{ if }x\in V_j
$$
with $\Phi^j_{n,\ell,q}$ given by (\ref{BBFunct2}). Taking Lemma \ref{LeMMaBBFund3}(i) into account (and recalling that $v$ is locally injective), it is easy to see that $\Psi_{n,\ell,q}$ is locally injective.  Using $(P_i)$ we can assert that
$$
I(\Psi_{n,\ell,q})\leq\int_{\Sigma}R_iW_0(\nabla\Psi_{n,\ell,q}(x))dx\hbox{ for all }n,\ell,q\geq 1.
$$
By Lemma \ref{LeMMaBBFund3}(ii) it is clear that for every $\ell,q\geq 1$, $\Psi_{n,l,q}\to v$ in $L^p(\Sigma;\RR^3)$. It follows that
$$
I(v)\leq \lim_{n\to+\infty}I(\Psi_{n,\ell,q})\leq\lim_{n\to+\infty}\int_{\Sigma}R_iW_0(\nabla\Psi_{n,\ell,q}(x))dx\hbox{ for all }\ell,q\geq 1.
$$
Moreover, from Lemma \ref{LeMMaBBFund3}(iii) we see that
$$
\lim_{q\to+\infty}\lim_{\ell\to+\infty}\lim_{n\to+\infty}\int_{\Sigma}R_iW_0(\nabla\Psi_{n,\ell,q}(x))dx=\int_\Sigma\R_{i+1}W_0(\nabla v(x))dx.
$$
Hence
$$
I(v)\leq \int_\Sigma\R_{i+1}W_0(\nabla v(x))dx,
$$
and the proof is complete.\hfill$\square$


\end{document}